\documentclass[11pt]{article}
\usepackage[utf8]{inputenc}
\usepackage{amsmath, amssymb}
\usepackage[margin=1in]{geometry}
\usepackage{hyperref}
\usepackage{enumitem}
\usepackage{parskip}
\usepackage{hyperref}

\usepackage{tabularx}
\usepackage{graphicx}
\usepackage{float}

\usepackage{booktabs,longtable,array}
\newcolumntype{P}[1]{>{\raggedright\arraybackslash}p{#1}}

\usepackage{booktabs,tabularx,array}
\newcolumntype{Y}{>{\raggedright\arraybackslash}X}

\usepackage{array}
\usepackage{float}
\renewcommand{\arraystretch}{1.4} 

\title{\textbf{Developing a Learner-Centered Teaching Routine
}}
\author{
Hyeeun Jang \\
Abilene Christian University \\
\texttt{hxj22bg@acu.edu}
}
\date{}  

\begin{document}

\maketitle

\begin{abstract}
This paper shares a classroom story from Fall 2022 to Spring 2025 about a learner centered routine in undergraduate mathematics. I use four steps: an opening question, a short mini lecture about meaning, structured small group work, and a short end of class exit check, sometimes with quick visuals. I used this mainly in elementary statistics, with some use in calculus and linear algebra. The evidence is local and practical: my notes, one minute exit checks, informal student comments and surveys, and conference feedback. Across courses, I saw less passive lecturing, more visible participation, and students explaining their reasoning during the closing time. Limits are clear: one instructor and no controlled comparisons. Next I plan to improve the prompts for different student groups and to simplify the visuals. I offer a simple routine and timing that other instructors can adapt to their own classes.
\end{abstract}

\section*{Introduction}

My students, especially non math majors, often struggle to connect with the subject. When I ask "Do you like math?" most say no unless they are majors. Many arrive with a negative view of mathematics. They think it is about memorizing formulas or solving complicated problems without knowing why, which makes learning feel distant and mechanical.I wanted to try something different. Instead of relying mainly on lectures and step by step solutions, I redesigned the course to give students more chances to talk, think, and explore ideas. My goal was to create a space where students could participate naturally, ask questions, and make sense of what they are learning. I also wanted them to feel ownership of their understanding, not just follow procedures. 

The approach is simple. Following Barbi (2019) and Lang (2016), I adopted brief focusing activities in the first five minutes to engage students and set the tone. These short routines activate prior knowledge, frame the day's work, and support retention. They are followed by a short mini lecture and small group discussion. I used this routine in an elementary statistics course, which often enrolls many non math majors, and I made small weekly adjustments toward a more learner centered approach. In my classes the routine had four parts: an opening focusing question, a brief meaning oriented mini lecture, structured small group work, and a brief whole class exit check, with occasional quick visuals when helpful.Small changes made a noticeable difference. Students were more engaged, more thoughtful, and more willing to share their reasoning. It was not only about getting the right answer; they explained why ideas worked and learned from one another. In what follows I describe how the routine looked in my courses and context, what I noticed from classroom notes and student comments, and a few next steps for this setting.

\section*{Learner-Centered Teaching Model}

This teaching routine consists of four parts that I used in my elementary statistics course: an opening question in the first five minutes, a brief meaning oriented mini lecture, structured small group work, and an end of class exit check. I did not change everything at once. Using these four elements, I made small adjustments each week, and I was surprised to see students begin to participate naturally and think for themselves. What follows comes from my classroom records in an elementary statistics course at a mid sized university, with about thirty to thirty five students per section, most of them non majors. Over time I kept the basic structure and applied it to other courses as well. In Linear Algebra and Calculus III, the same order (opening question, short mini lecture, and small group work) remained, with small adjustments for the content. For example, in Linear Algebra I paused to ask for common features across earlier examples so that students could build definitions together. In Calculus I, I checked one core idea each day and, when helpful, added quick visualizations in Mathematica or sometimes Python. In Statistics I used the exit check more often so that students explained not only how to compute but also why the result made sense. Short snapshots of these applications appear later in the section "Reflections and Outcomes".

\subsection*{1. First 5 minutes: Activating Curiosity}

Each class began with a short activity or question in the first five minutes. I wrote a simple but thought provoking open question on the board that connected to the day's topic. For example, before introducing box plots, I asked, "If a dataset has a mean of 100 and a median of 80, what might that look like?" That single question led some students to start talking right away, guessing ideas like skewed data or outliers before I said anything. In some classes, instead of writing on the board, I posted the question on Canvas, the learning management system at our university. I asked students to respond online within the first five minutes and then briefly read and comment on one or two classmates' posts. At times I also used a quick show of hands or invited students to write a single word on an index card to register a prediction before anyone spoke, which made it easier for quieter students to enter the conversation. This quick peer review style activity gave even quiet students a chance to engage thoughtfully and helped spark discussion before we moved into the mini lecture.

These brief activities helped students focus without pressure and shift into learning mode. They also encouraged curiosity and activated prior knowledge. I was inspired by Barbi (2019) and Lang (2016), who recommend using the beginning and end of class for structured thinking and engagement. In Spring 2023, when a few students still seemed unsure, I slowed the pacing and offered brief one to one check ins during office hours so they could join the opening question the next day. Later, in Calculus I, I used the same opening question as a short daily check, which gave me quick information for how to begin the next lesson.

The first five minutes often took a three part form: one concept in words, one short computation, and one brief interpretation of a small graph or situation. This balance kept both procedure and meaning in view and gave me quick information for grouping and pacing. In my notes, the share of students who completed this as intended rose from about 60\% in Fall 2023 to about 85\% in Fall 2024 as the routine became familiar.

\subsection*{2. Mini Lecture: Setting the Foundation}

After the opening question, I usually gave a traditional lecture for about 20 minutes. I did not try to cover every formula or term. I focused on key ideas that help students approach the concepts more meaningfully. For example, when I taught standard deviation, I emphasized what it means intuitively, how it shows the "typical distance" from the mean, and why it matters in real situations, rather than the exact computation. I kept the talk conversational and tied it to concrete examples. At times I paused to ask questions like "What do you notice?", "What would change if we doubled the sample?", and "Which part of this picture is fixed?" The goal was not just to transmit knowledge but to frame the problem so that the next stage would feel natural.

Sometimes I used simple analogies or everyday examples to make abstract ideas feel more accessible. For example, I explained sampling bias by comparing course ratings with online reviews and asking who tends to leave comments and who does not. Quick board sketches and very short visuals were enough to anchor the idea. When students seemed unsure, I slowed down, rephrased the explanation, or asked them to give their own examples. This approach clarified the content and helped students build their own models.

\begin{figure}[H]
  \centering
  \includegraphics[width=.5\linewidth]{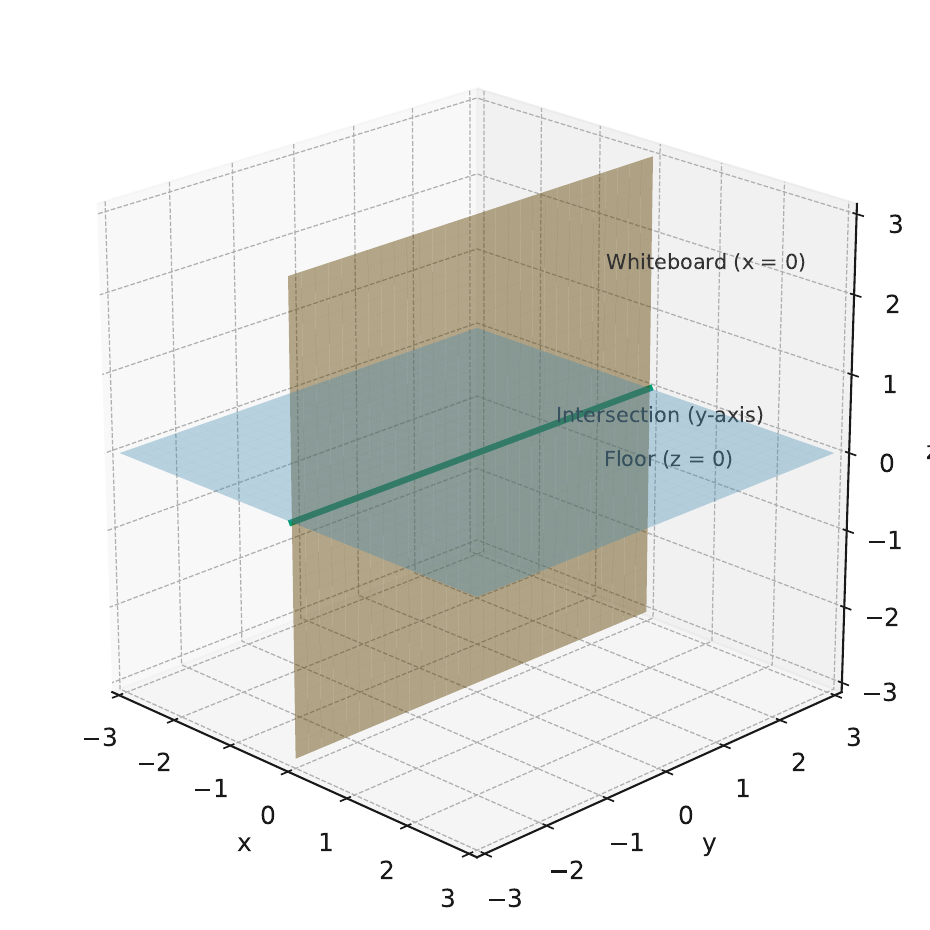}
  \caption{Two orthogonal planes $x=0$ (whiteboard) and $z=0$ (floor) intersecting along the $y$-axis.}
  \label{fig:whiteboard-floor}
\end{figure}

As the routine became familiar, I kept the basic timing mostly the same across courses (except in special cases). In Linear Algebra, to distinguish "being orthogonal" from "having only the zero vector in common," I showed a one minute Mathematica visualization of the "whiteboard" plane x=0 and the "floor" plane $z=0$ intersecting along the $y$-axis. This picture opened a discussion about subspace intersections and rank (see Figure~\ref{fig:whiteboard-floor}).

In Calculus I and III, I sometimes dragged a second point to show a secant slope approaching a tangent slope. Throughout, I kept the use of technology minimal and brief. I framed the idea, asked students to share and think on their own, and then moved to student discussion.

\subsection*{3. Group Discussion: Making Meaning Together}

After the mini lecture, students worked in pairs or small groups to solve problems or to discuss the concept learned that day. At first I set about twenty minutes for this stage. I walked among the groups, listened briefly, and asked light questions such as "What did you find here?" and "Did your group agree on a reason for that answer?" These short stops helped me see where students were thinking deeply and where they were stuck without putting anyone on the spot. I also provided a short worksheet with one or two guiding questions. At other times I said, "Try this example with your group and explain what is happening," and I encouraged students to lead their own group discussion actively. Early on, one or two strong voices sometimes led the group while quiet students only watched. To balance participation, I began rotating light roles such as reader, explainer, checker, and scribe, and I sometimes asked each student to give a short comment before the group reported. These small structures helped more students join the discussion.

As the semester went on, the mood of the room changed. Students asked more "why" questions, explained ideas to one another in their own words, and were willing to help when someone was stuck. Hesitant students found low pressure ways to take part, and their confidence grew. By the end of the term, I saw students not only finding answers but also comparing approaches, correcting each other gently, and noticing what they still wanted to understand. The change was subtle but steady. The class felt less like a performance and more like a shared exploration. This same group discussion routine carried over to my other courses with small adjustments. In Linear Algebra I used group time to surface common features across examples so that students could build definitions together. In Calculus I and III I paired short group work with a daily concept check so that students practiced both procedure and interpretation.

\subsection*{4. Last 5 Minutes: End of Class Check (Exit check) and Closure}

At the end of each class, I closed with a short end of class check (exit check) that lasted about four to five minutes. Students wrote one or two short responses about the activities from that day's class, for example one idea that became clear and one step their small group discussed. This was not a reflective essay but a simple check that gave students time to organize their thoughts and gave me a quick sense of what to revisit next time. When time was tight I used a show of hands or collected one sentence notes, and at other times I posted a single question on Canvas and reviewed a few responses before the next class.

Over time the exit check varied a little by course and semester. In Fall 2023 in elementary statistics I used one sentence notes on most days and took a few brief oral responses before the bell. In Spring 2024, while teaching two sections of statistics, I kept the same format but added a very short prompt on days with a new idea to confirm what was fixed and what was changing, and I used those notes to choose the next day's opening question. In Fall 2024 in Calculus I, I linked the closing to the opening question and, in the last minutes, had students confirm one key idea in words, check a small calculation, and add a brief interpretation in that order. When time was very limited I used a single yes or no question with a show of hands and continued in the next class. In Spring 2025 I connected the exit check to the day's short Mathematica visual by asking students to write one sentence about which picture from that day would help on the next problem.

Lang (2016) notes that a short end of class reflection supports memory and a sense of progress. I saw the same pattern in my courses. Students who were unsure during small group work often consolidated their understanding in the final minutes. In a quiet setting many used the time to organize the main point in their own words. I used the records from these exit checks to decide how to begin the next lesson. If several groups marked the same step as difficult, I opened with a short board example or a pair question. If many students named the same helpful idea, I used it as a bridge into the new task. Together with the opening question, the mini lecture, and the small group work, this closing helped keep a steady rhythm. Each part was simple, but the routine increased participation and brought more questions and explanations. For non majors, access and comfort improved, and interest and confidence rose overall.


\section*{Reflections and Outcomes}

Over the past several years I tried small adjustments, and paid attention to what helped my students in my courses and setting. The following is a brief account of that journey and what I learned in context. Unless noted otherwise, observations are \emph{approximate} and drawn from classroom notes, one-minute end-of-class exit checks, quick show-of-hands tallies, informal surveys, course evaluations, and simple logs. In practice I kept a brief record for each meeting (date, course, the small change I tried, and what I noticed). I saved photos of boards and collected exit cards to look for weekly patterns. I also kept quick tallies (e.g., number of hands raised for a claim; number of groups that reached an interpretation) so I could compare weeks without formal grading. 

\subsection*{Fall 2022: Early Missteps and Learning}
Fall 2022 was a humbling term. I assumed that detailed proofs and derivations would be enough. In Quantitative Reasoning and Elementary Statistics I relied on long lectures and traditional homework. Many nonmajors did not see the central goals, and evaluations often mentioned that the pace and level of abstraction were too high. I felt discouraged, but the feedback helped me see that accessibility, pacing, and guided practice mattered more than I had allowed. For the next term I began drafting clearer notes, planning short small-group prompts, and designing simpler study guides. Practically this meant adding a one-line main idea in the margin of each notes page, preparing two or three plain-language questions for each topic (e.g., what a large IQR tells us about spread), and creating study guides of five to seven items that paired one computation with one interpretation. I also scheduled office hours where students practiced one example before an exam.

\subsection*{Spring 2023: First Attempts, Mixed Success}
In Spring 2023 I introduced a short opening question and a brief problem-solving period in Quantitative Reasoning. When students explained ideas in their own words the tone of the room changed. Participation grew. At the same time my habits in Linear Algebra and Calculus kept those courses heavy and fast, and students still struggled with abstraction. I learned that working harder at the board would not be enough. I needed to slow down, invite student talk, and make space for interpretation as well as computation. Typical openers asked for a quick prediction or sketch. I used a visible 90-120 second timer and then had students pair-share before volunteers. About 5-10\% of students still had trouble joining group work, so I added brief one-on-one check-ins that often became an entry point to the next-day task.

\subsection*{Fall 2023: Beginning Group Discussion in Statistics}
In Fall 2023 I used short everyday prompts together with small-group work for about twenty minutes and then a whole-class share. Students explained to one another more freely and asked better questions. Written feedback often said that speaking first in groups and then hearing a short mini-lecture made understanding clearer. In my walk-around notes, hesitant students were more willing to contribute after a peer tried first. Some groups still waited for me to give the answer, but the room began to feel less like a performance and more like a space for making meaning together. I usually formed triads, gave light roles (reader, explainer, checker), and closed with a one-sentence claim. Prompts were often tiny datasets or everyday graphs with "what can we say with evidence?" I kept a simple table of common errors and addressed them in the next-day opening question.

\subsection*{Spring 2024: Greater Interaction, Persistent Challenges}
With two sections of Elementary Statistics I refined the structure: an \emph{opening question}, a brief \emph{mini-lecture} of about twenty minutes, \emph{small-group work} for about twenty minutes, and 5 minutes \emph{end-of-class exit check} that asked for a quick interpretation of a graph or everyday scenario. In Linear Algebra I emphasized dialogue around definitions. Before writing a new definition I asked what common features we noticed in earlier examples, and during proofs I paused to ask why a step was necessary. A more careful use of mathematical language began to take hold. Some topics (e.g., distinguishing the null space from the column space) remained difficult, and I noted to add visual support next term. Exit checks were brief: yes/no with a reason, a quick "what would change if..?"

\subsection*{Spring 2025: Bridging Abstraction with Visualization}
In Calculus and Linear Algebra I used short Mathematica and Python visuals to support intuition. For example, to distinguish "being orthogonal" from "sharing only the zero vector" I showed the planes $x=0$ (the "whiteboard") and $z=0$ (the "floor") intersecting along the $y$-axis. Students then named the subspaces and connected the picture to rank and intersections. Students said the graphics helped them link a symbolic statement to a concrete scene. In practice Mathematica let me produce quick classroom images with a few commands, and Python was helpful when I needed reproducibility. After each visual I asked for a single sentence that linked the picture to the symbols.

\subsection*{What Seemed to Change}
Shifts across terms are summarized in Table~\ref{tab:changes}. Group discussion became common and energetic, and many students turned in exit cards. Student initiated questions were more frequent and deeper than in 2022, and informal surveys suggested higher confidence, especially among nonmajors.

\small
\setlength{\tabcolsep}{4pt}
\renewcommand{\arraystretch}{1.15}

\begin{longtable}{>{\bfseries}P{4.1cm} P{2.9cm} P{2.9cm} P{2.9cm} P{2.9cm}}
\caption{Changes in Key Teaching Components}\label{tab:changes}\\
\toprule
Teaching Component & Fall 2022 & Fall 2023 & Fall 2024 & Spring 2025 \\
\midrule
\endfirsthead
\caption[]{Changes in Key Teaching Components (continued)}\\
\toprule
Teaching Component & Fall 2022 & Fall 2023 & Fall 2024 & Spring 2025 \\
\midrule
\endhead
\midrule
\multicolumn{5}{r}{Continued on next page}\\
\bottomrule
\endfoot
\bottomrule
\endlastfoot


Opening Question Engagement
& Not implemented
& Introduced, about 60\% active
& Consistent use, about 85\% active
& Routine; used as a daily check or to link a short visual, about 85\% to 90\% \\

Mini Lecture Focus
& Formula focused
& Shift to conceptual, 15 to 20 minutes
& Meaning driven and interactive
& Framed by quick context or visual; emphasis on interpretation \\

Small Group Participation
& Rare
& Encouraged, about 60\% active
& Regular and energetic, about 85\% or more
& Roles added in some weeks; broad participation sustained \\

Reflection / Exit Check Submission
& Not used
& 75\%
& 90\%
& 90\% (mostly one minute exit checks) \\

Critical Thinking in Reflections
& Not used
& Emerging (about 50\%)
& Rich responses (about 75\%)
& Interpretation with visuals and brief evidence based claims \\

Student Confidence (informal survey)
& Not measured
& Moderate (about 70\%)
& High (about 82\%)
& High (about 82\% to 85\%) \\

Peer Mentoring Observed
& Not observed
& Emerging
& Common and natural
& Common with occasional student led mini explanations \\

Student Initiated Questions
& Low
& Moderate increase
& Frequent and deeper
& Frequent with more why type questions and connections \\

Positive Course Feedback
& Standard
& Generally favorable
& Strongly positive
& Strongly positive with comments about clarity of visuals \\

Reported Math Anxiety
& High
& Some reduction
& Significant reduction
& Continued reduction (anecdotal) \\

\end{longtable}

\vspace{2pt}

 The learner-centered routine gave my students a voice and helped them feel like math was something they could understand and even enjoy. It also helped me become more reflective and intentional as a teacher. I still have more to learn, but this experience has shown me that even small changes can lead to a big difference.

\section*{Limitations and Future Directions}

 In Fall 2022 the approach was mostly lecture based. In Fall 2023 and Fall 2024 I gradually introduced the opener, the short mini lecture, small group work, and the closing exit check, and in Spring 2025 I added brief visuals. Because the routine matured step by step, it is difficult to separate which element had the greatest impact. Section size, time of day, and the mix of majors and nonmajors also varied across semesters and may have shaped the outcomes. Some students preferred individual study and a faster pace, so a group based routine was not equally beneficial for everyone. My records focused mainly on participation, interpretation, and classroom climate. Visuals were used deliberately in a minimal way, which worked in my setting but may not be appropriate for every context. The percentages in Table~\ref{tab:changes} are approximate and should be read as descriptive classroom notes rather than generalizable findings. Even so, the changes I observed were encouraging. Students were more willing to take small intellectual risks, to talk with one another, and to explain ideas in their own words. Especially among nonmajors, there were informal reports of reduced anxiety and clearer sense of meaning in results. These signs suggest that the routine is worth exploring in a more systematic way.

\subsection*{Future Directions}

In the near future I plan to keep the routine simple while continuing to keep records and to add a few light tools that make comparisons across terms possible.

\begin{itemize}
\item Continue to track simple counts such as opener participation rate, the proportion of groups that reach an interpretation, and the number of student initiated questions, while protecting privacy.
\item Work with colleagues to try the routine in other entry level courses such as college algebra or precalculus, and note what parts need adjustment.
\item Build a small open collection of materials that match the routine, including a page of opener prompts, a page of exit check stems, and a few simple Mathematica or Python notebooks for visuals that students found helpful.
\item Continue small equity moves such as quiet think time before discussion and rotating who speaks first, and check whether these broaden participation.
\item Explore simple digital tools that support quick visuals and formative checks without taking over class time, and record both benefits and tradeoffs.
\end{itemize}

In the longer term, if these steps prove useful, I hope to collaborate with colleagues on a study that includes common tasks across sections, short affect measures of confidence, and classroom observations focused on interpretation and explanation. Even in that work I intend to keep the emphasis on small and concrete routines that fit the daily realities of entry level courses. The goal is not only to know whether students can compute, but also when and why they feel ready to explain a result in context.

\section*{References}

\begin{itemize}
    \item Lang, J. M. (2016). \emph{Small Teaching: Everyday Lessons from the Science of Learning}. John Wiley \& Sons.
\item Barbi, H. (2019). ``Three focusing activities to engage students in the first five minutes of class." \emph{Faculty Focus Higher Ed Teaching \& Learning}. \href{https://www.facultyfocus.com/articles/blended-flipped-learning/three-focusing-activities-engage-students-first-five-minutes-class}{https://www.facultyfocus.com/articles/blended-flipped-learning/three-focusing-activities-engage-students-first-five-minutes-class}

\end{itemize}

\end{document}